\documentclass[11pt, reqno, psamsfonts]{amsart}
\pdfoutput=1

\usepackage{amssymb}
\usepackage{amsthm}
\usepackage{amsmath}
\usepackage{latexsym}
\usepackage[T1]{fontenc}
\usepackage[utf8]{inputenc}
\usepackage[russian, french, english]{babel}

\usepackage{graphicx}
\usepackage{wrapfig}
\usepackage[justification=centering, labelfont=bf]{caption}
\usepackage{mathtools}
\usepackage{amsbsy}
\usepackage[inline]{enumitem}
\usepackage{mathrsfs}
\usepackage{array}
\usepackage{multicol}
\usepackage{stmaryrd}
\usepackage{cancel}
\usepackage{lmodern}
\usepackage{mathabx}
\usepackage{upgreek}
\usepackage{titlesec}
\usepackage{titletoc}
\usepackage[spacing=true,kerning=true,babel=true,tracking=true]{microtype}
\usepackage[shortcuts]{extdash}
\usepackage[foot]{amsaddr}
\usepackage[left=1in,right=1in,top=1in,bottom=1in,bindingoffset=0cm]{geometry}
\usepackage{bm}
\usepackage{centernot}
\usepackage{mdframed}
\usepackage[hidelinks]{hyperref}
\usepackage{xspace}
\usepackage[skins]{tcolorbox}
\usepackage{framed}
\usepackage[justification=centering, labelfont=bf]{caption}
\usepackage{tikz}
\usetikzlibrary{shapes,snakes}
\usetikzlibrary{arrows.meta}
\usetikzlibrary{decorations.pathmorphing}
\usetikzlibrary{patterns}
\usepackage{float}

\usepackage[
    sortcites,
    backend=biber, style=alphabetic, sorting=nyt, maxnames=100,backref=true]{biblatex}
\title{\sffamily Embedding Borel graphs into grids of\\ asymptotically optimal dimension}
\date{}

\author{Anton~Bernshteyn}
\address{\normalfont (AB) Department of Mathematics, University of California, Los Angeles, CA, USA}
\email{bernshteyn@math.ucla.edu}

\author{Jing~Yu}
\address{\normalfont (JY) Shanghai Center for Mathematical Sciences, Fudan University, Shanghai, China}
\email{jyu@fudan.edu.cn}

\thanks{AB's research is partially supported by the NSF CAREER grant DMS-2528522 and the Sloan Research Fellowship (2025). JY's research is partially supported by the China Postdoctoral Science Foundation grant GZC20252041 and the National Natural Science Foundation of China grant 12371343 (PI: Hehui Wu).}

\newtheoremstyle{bfnote}%
{}{}%
{\slshape}{}%
{\bfseries}{\bfseries.}%
{ }%
{\thmname{#1}\thmnumber{ #2}\thmnote{ \ep{\normalfont{}#3}}}

\theoremstyle{bfnote}
\newtheorem{theo}{Theorem}[section]
\newtheorem*{theo*}{Theorem}
\newtheorem{prop}[theo]{Proposition}
\newtheorem{lemma}[theo]{Lemma}

\newtheorem{corl}[theo]{Corollary}

\newtheorem*{corl*}{Corollary}

\theoremstyle{definition}
\newtheorem{defn}[theo]{Definition}
\newtheorem*{defn*}{Definition}
\newtheorem{exmp}[theo]{Example}

\newtheorem*{exmp*}{Example}

\theoremstyle{remark}
\newtheorem*{ques*}{Question}
\newtheorem*{remk*}{Remark}

\newcommand*{\myproofname}{Proof}

\makeatletter
\newcommand{\neutralize}[1]{\expandafter\let\csname c@#1\endcsname\count@}
\makeatother

 \newenvironment{lemmacopy}[1]
{\renewcommand{\thetheo}{\ref*{#1}}%
	\neutralize{theo}\phantomsection
	\begin{lemma}}
	{\end{lemma}}

\newcommand{\set}[1]{\{#1\}}
\newcommand{\N}{{\mathbb{N}}}
\newcommand{\Z}{\mathbb{Z}}

\newcommand{\R}{\mathbb{R}}

\renewcommand{\epsilon}{\varepsilon}

\renewcommand{\phi}{\varphi}
\renewcommand{\theta}{\vartheta}
\renewcommand{\leq}{\leqslant}
\renewcommand{\geq}{\geqslant}
\newcommand{\defeq}{\coloneqq}

\newcommand{\bemph}[1]{{\normalfont#1}} 
\newcommand{\ep}[1]{\bemph{(}#1\bemph{)}} 

\newcommand{\emphd}[1]{{\fontseries{b}\selectfont\textsf{#1}}}
\newcommand{\acts}{\mathrel{\reflectbox{$\righttoleftarrow$}}}

\newcommand{\G}{\Gamma}
\newcommand{\pto}{\dashrightarrow}
\newcommand{\dom}{\mathsf{dom}}
\newcommand{\gr}[1]{\Z^{#1}_\infty}
\newcommand{\shgr}{\mathsf{ShiftGrid}}
\newcommand{\Free}{\mathsf{Free}}

\newcommand{\asi}{\mathsf{asi}}

\newcommand{\asdimB}{\mathsf{asdim}_\mathsf{B}}
\newcommand{\embB}{\emb_\mathsf{B}}
\newcommand{\dist}{\mathsf{dist}}

\newcommand{\emb}{\mathsf{emb}}
\newcommand{\Sch}{\mathsf{Sch}}
\renewcommand{\deg}{\mathsf{deg}}

\numberwithin{equation}{section}

\newenvironment{scproof}[1][]{\begin{proof}[\textsc{\upshape{Proof}}#1]}{\end{proof}}

\titleformat{\section}[block]{\large\bfseries\sffamily}{\thesection.}{1ex}{}
\titleformat{\subsection}[block]{\bfseries\sffamily}{\thesubsection.}{1ex}{}
\titleformat{\subsubsection}[block]{\itshape}{\bfseries\upshape\sffamily\thesubsubsection.}{1ex}{}

\titlespacing*{\section}{0pt}{*3}{*1}
\titlespacing*{\subsection}{0pt}{*3}{*1}
\titlespacing*{\subsubsection}{0pt}{*2}{*1}

\titlecontents{section}
 [1.5em] %
 {\smallskip}
 {\bfseries\thecontentslabel\hspace{1.02em}}
{\bfseries}
 {\,\,\titlerule*[0.77pc]{}\bfseries\contentspage}
\titlecontents{subsection}
 [4em] %
 {\smallskip}
 {\thecontentslabel\hspace{1.02em}}
{\hspace*{2.32em}}
 {\,\,\titlerule*[0.77pc]{.}\contentspage}

\renewbibmacro{in:}{}

\renewbibmacro*{volume+number+eid}{%
	\printfield{volume}%
	\setunit*{\addnbspace}
	\printfield{number}%
	\setunit{\addcomma\space}%
	\printfield{eid}}

\DeclareFieldFormat[article]{volume}{\textbf{#1}\space}
\DeclareFieldFormat[article]{number}{\mkbibparens{#1}}

\DeclareFieldFormat{journaltitle}{#1,}
\DeclareFieldFormat[thesis]{title}{\mkbibemph{#1}\addperiod}
\DeclareFieldFormat[article, unpublished, thesis]{title}{\mkbibemph{#1},}
\DeclareFieldFormat[book]{title}{\mkbibemph{#1}\addperiod}
\DeclareFieldFormat[unpublished]{howpublished}{#1, }

\DeclareFieldFormat{pages}{#1}

\DeclareFieldFormat[article]{series}{Ser.~#1\addcomma}

\setlist{topsep=3pt,itemsep=3pt}

\begin{document}


    \maketitle
    
    
    \begin{abstract}
        Let $G$ be a Borel graph all of whose finite subgraphs embed into the $d$\-/dimensional grid with diagonals. We show that then $G$ itself admits a Borel embedding into the Schreier graph of a free Borel action of $\Z^{O(d)}$. This strengthens an earlier result of the authors, in which $O(d)$ is replaced by $O(\rho \log \rho)$, where $\rho$ is the polynomial growth rate of $G$.
    \end{abstract}
    
    \section{Introduction}

    In this paper we continue our study of the geometry of Borel graphs of polynomial growth initiated in \cite{BernshteynYu}. All graphs in this paper are undirected and simple; that is, a graph $G$ consists of a set of vertices $V(G)$ and a set of edges $E(G) \subseteq [V(G)]^2$, where $[X]^2$ is the set of all $2$-element subsets of $X$. We use the word ``metric'' to mean ``extended metric,'' i.e., we allow distances in a metric space to be infinite. For a graph $G$, we let $\dist_G$ be the \emphd{graph metric} on $V(G)$, where $\dist_G(u,v)$ is the minimum number of edges in a $uv$-path in $G$ if such a path exists and $\infty$ otherwise. We also let $B_G(u,r)$ be the closed ball of radius $r$ around $u \in V(G)$ in the graph metric. A graph $G$ is \emphd{of polynomial growth} if there is $\rho \in [0,\infty)$ such that $|B_G(u,r)| \leq (r+1)^\rho$ for all $u \in V(G)$ and $r \geq 1$ (in particular, every graph of polynomial growth is locally finite). We call the smallest such $\rho$ the \emphd{growth rate} of $G$ and denote it by $\rho(G)$.\footnote{In \cite{BernshteynYu} this parameter is called the \emph{exact} growth rate of $G$, to be contrasted with its \emph{asymptotic} growth rate.} An \emphd{embedding} of a graph $G$ into a graph $H$ is an injective mapping $f \colon V(G) \to V(H)$ such that $\set{f(u),f(v)} \in E(H)$ for all $\set{u,v} \in E(G)$. 

    \begin{wrapfigure}[10]{r}{0.3\textwidth}
			\centering
			\begin{tikzpicture}
				
				\begin{scope}[scale=1.6]
				    \filldraw (0,0) circle (2pt);
				    \filldraw (0.5,0) circle (2pt);
				    \filldraw (1,0) circle (2pt);
				    \filldraw (1.5,0) circle (2pt);
				    \filldraw (0,0.5) circle (2pt);
				    \filldraw (0.5,0.5) circle (2pt);
				    \filldraw (1,0.5) circle (2pt);
				    \filldraw (1.5,0.5) circle (2pt);
				    \filldraw (0,1) circle (2pt);
				    \filldraw (0.5,1) circle (2pt);
				    \filldraw (1,1) circle (2pt);
				    \filldraw (1.5,1) circle (2pt);
				    \filldraw (0,1.5) circle (2pt);
				    \filldraw (0.5,1.5) circle (2pt);
				    \filldraw (1,1.5) circle (2pt);
				    \filldraw (1.5,1.5) circle (2pt);
				    
				    \draw (-0.2,0) -- (1.7,0) (-0.2,0.5) -- (1.7,0.5) (-0.2,1) -- (1.7,1) (-0.2,1.5) -- (1.7,1.5) (0,-0.2) -- (0,1.7) (0.5,-0.2) -- (0.5, 1.7) (1,-0.2) -- (1,1.7) (1.5,-0.2) -- (1.5, 1.7);
				    
				    \draw (-0.2,-0.2) -- (1.7,1.7) (-0.2, 0.3) -- (1.2,1.7) (-0.2, 0.8) -- (0.7,1.7) (-0.2, 1.3) -- (0.2, 1.7) (0.3,-0.2) -- (1.7,1.2) (0.8,-0.2) -- (1.7,0.7) (1.3,-0.2) -- (1.7,0.2);
				    
				    \draw (1.5+0.2,-0.2) -- (1.5-1.7,1.7) (1.5+0.2, 0.3) -- (1.5-1.2,1.7) (1.5+0.2, 0.8) -- (1.5-0.7,1.7) (1.5+0.2, 1.3) -- (1.5-0.2, 1.7) (1.5-0.3,-0.2) -- (1.5-1.7,1.2) (1.5-0.8,-0.2) -- (1.5-1.7,0.7) (1.5-1.3,-0.2) -- (1.5-1.7,0.2);
				\end{scope}
			\end{tikzpicture}
			\caption{A fragment of $\gr{2}$.\label{fig:grids}}
	\end{wrapfigure}

    All finite graphs are of polynomial growth, although their growth rates may be arbitrarily high. A more interesting example is the \emphd{$d$-dimensional grid with diagonals} 
    $\gr{d}$,\footnote{This notation is taken from \cite{linial1995geometry,krauthgamer2007intrinsic}. In \cite{BernshteynYu}, the notation $\mathsf{Grid}_{d,\infty}$ is used instead.} which is defined by
    \[\displaywidth=\parshapelength\numexpr\prevgraf+2\relax
        \begin{array}{rl}
        V(\gr{d}) &\defeq\, \Z^d,\\
        E(\gr{d}) &\defeq\, \set{\set{u,v} \in [\Z^d]^2 \,:\, \|u - v\|_\infty = 1}.
        \end{array}
    \]
    That is, $\gr{d}$ is the Cayley graph of $\Z^d$ associated to the generating set $S_{d,\infty} \defeq \set{k \in \Z^d : \|k\|_\infty = 1}$. Note that $\rho(\gr{d}) = \Theta(d)$.
    
    Since grids with diagonals contain arbitrarily large cliques, every finite graph $F$ embeds into $\gr{d}$ for some $d \in \N$.  
    Answering a question of Levin and Linial, London, and Rabinovich \cite{linial1995geometry}, Krauthgamer and Lee \cite{krauthgamer2007intrinsic} showed that in fact $F$ embeds into a grid with diagonals whose dimension depends only on $\rho(F)$: 

    \begin{theo}[{Krauthgamer--Lee \cite[Theorem 5.8]{krauthgamer2007intrinsic}}]\label{theo:KL1}
        If $F$ is a finite graph of growth rate $\rho$, then $F$ has an embedding into $\gr{d}$ for some $d = O(\rho \log \rho)$. 
    \end{theo}


    Theorem~\ref{theo:KL1} motivates the following definition:

    \begin{defn}[Embedding dimension]
        The \emphd{embedding dimension}\footnote{In \cite{BernshteynYu} we used the term ``injective dimension'' instead.} of a finite graph $F$, in symbols $\emb(F)$, is the least $d \in \N$ such that $F$ admits an embedding into $\gr{d}$. 
        For an infinite graph $G$, we let $\emb(G)$ be the supremum of $\emb(F)$ taken over all finite subgraphs $F \subseteq G$.
    \end{defn}

    Combining the relation $\rho(\gr{d}) = \Theta(d)$ with Theorem~\ref{theo:KL1}, we obtain the following asymptotic bounds on $\emb(G)$ in terms of $\rho(G)$ for any graph $G$ of polynomial growth:
    \begin{equation}\label{eq:bounds}
        \rho(G) \,\preccurlyeq\, \emb(G) \,\preccurlyeq\, \rho(G) \log\rho(G).
    \end{equation}
    (Here ``$X \preccurlyeq Y$'' stands for ``$X = O(Y)$.'') Levin and Linial, London, and Rabinovich conjectured that the lower bound in \eqref{eq:bounds} is always tight, i.e., $\emb(G) = \Theta(\rho(G))$ for every graph $G$ of polynomial growth \cite[Conjecture 8.2]{linial1995geometry}. This was disproved by Krauthgamer and Lee, who used expander graphs to construct examples with $\emb(G) = \Theta(\rho(G) \log \rho(G))$, which shows that the upper bound in \eqref{eq:bounds} cannot be improved in general 
    \cite[Theorem 3.2]{krauthgamer2007intrinsic}. Nevertheless, for certain classes of graphs, the Levin/Linial--London--Rabinovich conjecture does hold. For example, the following results were established by Krauthgamer and Lee: 


    \begin{theo}[{Krauthgamer--Lee \cite[Corollary 2.10 and Theorem 4.10]{krauthgamer2007intrinsic}}]\label{theo:KL2}
        Let $G$ be a graph of polynomial growth with growth rate $\rho$. Fix integers $\ell \geq 3$ and $s \geq 1$.
        
        \begin{enumerate}[label=\ep{\normalfont{\arabic*}}]
            \item\label{item:cycles} If $G$ has no induced cycle of length at least $\ell$, then $\emb(G) = O((\log \ell) \, \rho)$.

            \item\label{item:minors} If $G$ has no minor isomorphic to the complete bipartite graph $K_{s,s}$, then $\emb(G) = O(4^{s} \rho)$.
        \end{enumerate}
    \end{theo}

    Note that since every finite graph $H$ is a minor of $K_{s,s}$ for some $s \in \N$,  Theorem~\ref{theo:KL2}\ref{item:minors} implies that $\emb(G) = O_H(\rho(G))$ whenever $G$ is $H$-minor-free.

    We are interested in extending the above results to the realm of Borel graphs. The systematic study of combinatorial properties of Borel graphs was launched by Kechris, Solecki, and Todorcevic in their seminal paper \cite{KST} and has since developed into a rich subject with many connections to other areas such as dynamical systems and computer science; see \cite{KechrisMarks} for a survey by Kechris and Marks and \cite{Pikh_survey} for an introductory article by Pikhurko. We also direct the reader to \cite{KechrisDST,AnushDST} for general descriptive set theory background. A useful informal perspective on Borel combinatorics (at least in the locally finite regime) is that instead of working with a single connected graph, we are given a family of graphs, and our goal is to solve a combinatorial problem on all of them in a ``uniform'' fashion.

    \begin{exmp}[Schreier graphs]
        Many examples in Borel combinatorics come from actions of finitely generated groups. Let $\G = \langle S \rangle$ be a group generated by a finite symmetric set $S$ that does not contain the identity element and suppose that $\G \acts X$ is an action of $\G$ on a standard Borel space $X$ by Borel automorphisms.\footnote{The reader unfamiliar with standard Borel spaces may assume that $X$ is, say, the unit interval $[0,1]$ or the real line $\R$. This usually results in no loss of generality, thanks to the Borel isomorphism theorem \cite[Theorem 15.6]{KechrisDST}.} Assume additionally that the action $\G \acts X$ is \emphd{free}, meaning that the stabilizer of every point $x \in X$ is trivial. Then we can form the \emphd{$S$-Schreier graph} $\Sch(X, S)$ of this action by pasting a copy of the Cayley graph of $\G$ onto each orbit. More precisely, we let
        \[
            \begin{array}{rl}
        V(\Sch(X,S)) &\defeq\, X,\\
        E(\Sch(X,S)) &\defeq\, \set{\set{x,\sigma \cdot x} \,:\, x \in X, \, \sigma \in S}.
        \end{array}
        \]
        The fact that the action $\G \acts X$ is free implies that each component of $\Sch(X,S)$ is isomorphic to the Cayley graph of $\G$ associated to the generating set $S$. Thus, we can think of $\Sch(X,S)$ as a (typically uncountable) family of copies of the Cayley graph ``put together'' in a single space.
    \end{exmp}

    The general definition is as follows:

    \begin{defn}[Borel graphs]\label{defn:Borel_graph}
        A \emphd{Borel graph} is a graph $G$ whose vertex set $V(G)$ is a standard Borel space and whose edge set $E(G)$ is a Borel subset of $[V(G)]^2$.
    \end{defn}

    It is not hard to see that Theorem~\ref{theo:KL1} can be extended from finite graphs $F$ to arbitrary connected graphs of polynomial growth. It follows that if $G$ is a Borel graph of polynomial growth, then every component of $G$ embeds into a grid with diagonals. Can this embedding be achieved for all components of $G$ in a ``uniform'' way? More precisely, we seek a Borel embedding of $G$ into a Borel graph whose components are copies of $\gr{d}$. 
    Recall that $S_{d,\infty} = \set{k \in \Z^d \,:\, \|k\|_\infty = 1}$.

    \begin{defn}[Borel embedding dimension]\label{defn:Borel_dim}
        We say that a Borel graph $H$ is a \emphd{Borel $\gr{d}$-graph} if it is the $S_{d,\infty}$-Schreier graph of some free Borel action $\Z^d \acts V(H)$. The \emphd{Borel embedding dimension} of a Borel graph $G$, in symbols $\embB(G)$, is the least $d \in \N$ such that $G$ has a Borel embedding 
        into a Borel $\gr{d}$-graph. 
        If no such $d \in \N$ exists, we let $\embB(G) \defeq \infty$.
    \end{defn}

    It is possible to alter Definition \ref{defn:Borel_dim} by considering \emph{strong} embeddings, i.e., those in which vertices in different components of $G$ are mapped to distinct components of $H$; this would not affect the required value of $d$---see Proposition~\ref{prop:strong}. Also, instead of working with arbitrary Borel $\gr{d}$-graphs, one may only focus on certain ``universal'' examples, as explained in Proposition~\ref{prop:universal}. 
    
    It is clear from the definition that $\embB(G) \geq \emb(G)$. One of the central results of \cite{BernshteynYu} is that the asymptotic upper bound in Theorem~\ref{theo:KL1} holds for the Borel embedding dimension:

    \begin{theo}[{AB--JY \cite[Theorem 1.10]{BernshteynYu}}]\label{theo:BY}
        If $G$ is a Borel graph of polynomial growth with growth rate $\rho$, then $\embB(G) = O(\rho \log \rho)$.
    \end{theo}

    Of course, if $\emb(G) = \Theta(\rho \log \rho)$ (for instance when $G$ has one of the examples from \cite[Theorem 3.2]{krauthgamer2007intrinsic} as a subgraph), the bound in Theorem~\ref{theo:BY} cannot be improved. On the other hand, as exemplified by Theorem~\ref{theo:KL2}, there exist interesting classes of graphs with $\emb(G) = o(\rho \log \rho)$, and a natural question is whether it is possible to achieve an asymptotically stronger upper bound on $\embB(G)$ for graphs in such classes. Our main result is that the answer is ``yes''; in fact, we show that the Borel embedding dimension of a Borel graph $G$ never exceeds its ordinary embedding dimension by more than a constant factor:
    
    \begin{theo}[Main result: Borel embedding dimension]\label{theo:main}
        There exists an absolute constant $C > 0$ such that if $G$ is a Borel graph of polynomial growth, then $\embB(G) \leq C\,\emb(G)$.
    \end{theo}

    For example, combining Theorem~\ref{theo:main} with Theorem~\ref{theo:KL2} yields the following bounds: 

    \begin{corl}[Induced cycles and forbidden minors]\label{corl:minors}
        Let $G$ be a Borel graph of polynomial growth with growth rate $\rho$. Fix integers $\ell \geq 3$ and $s \geq 1$.
        
        \begin{enumerate}[label=\ep{\normalfont{\arabic*}}]
            \item\label{item:cycles} If $G$ has no induced cycle of length at least $\ell$, then $\embB(G) = O((\log \ell) \, \rho)$.

            \item\label{item:minors} If $G$ has no minor isomorphic to $K_{s,s}$, then $\embB(G) = O(4^{s} \rho)$.
        \end{enumerate}
    \end{corl}

    A version of Corollary~\ref{corl:minors}\ref{item:minors} with no explicit dependence on $s$ appeared earlier in the second named author's PhD thesis with a different proof \cite{Jing}.

    An important consideration in the study of metric embeddings is the degree to which they affect the distances between points. In that direction, we prove the following stronger form of Theorem~\ref{theo:main} that yields embeddings that do not reduce distances between vertices ``too much'': 

    \begin{theo}[Embeddings with low distance reduction]\label{theo:coarse}
        Let $G$ be a Borel graph of polynomial growth with $\emb(G) = d$. For each $\epsilon > 0$, there exists a Borel embedding $f \colon V(G) \to V(H)$ of $G$ into a Borel $\gr{n}$-graph $H$ such that $n \leq C_\epsilon\, d$ and for all $u$, $v \in V(G)$, 
        \begin{equation}\label{eq:lower_bound_on_distance}
            \dist_G(u,v) \,\geq\, R_{\epsilon,d} \quad \Longrightarrow \quad \dist_H(f(u), f(v)) \,\geq\, \dist_G(u,v)^{1 - \epsilon}.
        \end{equation} Here $C_\epsilon > 0$ depends only on $\epsilon$ and $R_{\epsilon,d} > 0$ depends only on $\epsilon$ and $d$.
    \end{theo}
    
    To clarify, implication \eqref{eq:lower_bound_on_distance} in Theorem~\ref{theo:coarse} is satisfied even if $\dist_G(u,v) = \infty$, in which case it asserts that $\dist_H(f(u),f(v)) = \infty$ as well. Theorem~\ref{theo:main} trivially follows from Theorem~\ref{theo:coarse} by fixing an arbitrary constant $\epsilon > 0$ and letting $C \defeq C_\epsilon$.

    We made no attempt to optimize the constant $C$ in Theorem~\ref{theo:main}. Our argument gives $C \approx 10^7$, but this can almost certainly be improved with greater care in the calculations. However, that would require reworking some of the rather technical analysis in \cite{BernshteynYu}, used as a black box in this paper, and at any rate it is unlikely that this method can yield an optimal bound. 
    Note that there do exist Borel graphs $G$ of polynomial growth with $\embB(G) > \emb(G)$; for example, an undirectable Borel forest of lines $L$ satisfies $\emb(L) = 1 < \embB(L)$ (see \cites[Remark 6.8]{KechrisMiller}{MillerThesis}).

    To conclude the introduction, let us say a few words about the proof of Theorem~\ref{theo:coarse}. In addition to the results and techniques from \cite{BernshteynYu}, it relies on one key new lemma. 
    A function $f \colon X \to Y$ between metric spaces $(X, d_X)$ and $(Y, d_Y)$ is
    \begin{itemize}
        \item \emphd{$k$-Lipschitz} if $d_Y(f(x), f(y)) \leq k\,d_X(x,y)$ for all $x$, $y \in X$,
        \item \emphd{$R$-locally injective} if $f(x) \neq f(y)$ for all $x$, $y \in X$ such that $0 < d_X(x,y) \leq R$.
    \end{itemize}
    The two main types of metric spaces we consider in this paper are $(\Z^d, \|\cdot\|_\infty)$ and $(V(G), \dist_G)$ for a graph $G$; these are the spaces we mean when the metric is not specified explicitly.

    \begin{lemma}[Key lemma]\label{lemma:key}
        Let $G$ be a Borel graph of polynomial growth with $\emb(G) = d$. Then, for any $R > 0$, there exists an $R$-locally injective $1$-Lipschitz Borel function $f \colon V(G) \to \Z^{4d}$.
    \end{lemma}

    In our construction of a Borel embedding witnessing the bound on $\embB(G)$, we use Lemma~\ref{lemma:key} to ensure injectivity for vertices that are close to each other, while vertices that are far apart are handled using our earlier result, namely \cite[Theorem 5.7]{BernshteynYu}. The details of the construction are presented in \S\ref{sec:merge}. The proof of Lemma~\ref{lemma:key} is given in \S\ref{sec:asi}. An important tool used there is the asymptotic separation index of Borel graphs, which is a powerful concept that was introduced by Conley, Jackson, Marks, Seward, and Tucker-Drob \cite{Dimension} and has found numerous applications in Borel combinatorics \cite{BWKonig,Dimension,ASIalgorithms,FelixFinDim,BernshteynWeilacher}. We finish the paper with some further remarks concerning embeddings into Borel $\gr{d}$-graphs in \S\ref{sec:further}.

    \section{Proof of Theorem~\ref{theo:coarse} assuming Lemma~\ref{lemma:key}}\label{sec:merge}

    In this section we prove Theorem~\ref{theo:coarse} (and hence also Theorem~\ref{theo:main}) assuming the key Lemma~\ref{lemma:key}. We shall use the following notation for actions of $\Z^n$: if $\Z^n \acts X$ is an action, we write $k \bm{+} x$ to mean the result of acting on a point $x \in X$ by a group element $k \in \Z^n$. In addition to Lemma~\ref{lemma:key}, we need the following fact, obtained by applying \cite[Theorem 5.7]{BernshteynYu} with $b = \rho$ and $r = 1$:

    \begin{theo}[{AB--JY \cite[Theorem 5.7]{BernshteynYu}}]\label{theo:coarse_old}
        Let $G$ be a Borel graph of polynomial growth with growth rate at most $\rho \geq 1$. For each $0 < \epsilon < 1/2$, there is a $1$-Lipschitz Borel map $f \colon V(G) \to V(H)$ from $G$ to a Borel $\gr{n}$-graph $H$ such that $n = \left\lceil 10^7  \log(1/\epsilon)\epsilon^{-2} \, \rho\right\rceil$ and for all $u$, $v \in V(G)$,
        \begin{equation}\label{eq:far}
            \dist_G(u,v) \,\geq\, \left(10^7 \epsilon^{-1}\, \rho\right)^{\frac{3000}{\epsilon^2}} \quad \Longrightarrow \quad \dist_H(f(u), f(v)) \,\geq\, \dist_G(u, v)^{1-\epsilon}.
        \end{equation}
    \end{theo}




        Now let $G$ be a Borel graph of polynomial growth with $\emb(G) = d$ and fix $\epsilon > 0$. 
        Without loss of generality, we may assume that $d \geq 1$ and $\epsilon < 1/2$. Set $\rho \defeq \max\set{\rho(G), 1}$ and 
        \[
            R \,\defeq\, \left(10^7 \epsilon^{-1}\, \rho\right)^{\frac{3000}{\epsilon^2}},
        \]
        and let $n$, $H$, and $f$ be given by Theorem~\ref{theo:coarse_old}. 
        Also, let $\Z^n \acts V(H)$ be the free Borel action that generates $H$. Note that $f$ \emph{almost} satisfies the requirements of Theorem~\ref{theo:coarse}; the only property $f$ may lack is injectivity: although $\dist_H(f(u),f(v))$ is large when $u$ and $v$ are far apart, it is possible that $f(u) = f(v)$ when $\dist_G(u,v) < R$. This is precisely the issue addressed by Lemma~\ref{lemma:key}, which yields an $R$-locally injective $1$-Lipschitz Borel function $h \colon V(G) \to \Z^{4d}$. To combine $f$ and $h$, we define an action $\Z^{n + 4d} \acts V(H) \times \Z^{4d}$ by the formula
        \[
            (k, \ell) \bm{+} (x,m) \,\defeq\, (k \bm{+} x, \ell + m) \quad \text{for all $x \in V(H)$, $m \in \Z^{4d}$, $k \in \Z^n$, and $\ell \in \Z^{4d}$}.
        \]
        Since both the action $\Z^n \acts V(H)$ and the translation action of $\Z^{4d}$ on itself are free, we conclude that the action $\Z^{n + 4d} \acts V(H) \times \Z^{4d}$ is free as well. Let $H^*$ be the Borel $\gr{n + 4d}$-graph generated by this action. It is easy to see that for all $(x, m)$, $(x', m') \in V(H) \times \Z^{4d} = V(H^*)$,
        \begin{equation}\label{eq:infty_norm}
            \dist_{H^*}((x,m),(x',m')) \,=\, \max\set{\dist_H(x,x'),\, \|m - m'\|_\infty}.
        \end{equation}
        Define $f^* \colon V(G) \to V(H^*) = V(H) \times \Z^{4d}$ by
        \[
            f^*(u) \,\defeq\, (f(u), h(u)).
        \]
        Then, by \eqref{eq:infty_norm}, for any two vertices $u$, $v \in V(G)$,
        \[
            \dist_{H^*}(f^*(u), f^*(v)) \,=\, \max\set{\dist_H(f(u),f(v)),\, \|h(u) - h(v)\|_\infty}.
        \]
        Since both $f$ and $h$ are $1$-Lipschitz, it follows that $f^*$ is $1$-Lipschitz as well. Now consider any $u$, $v \in V(G)$. If $0 < \dist_G(u,v) \leq R$, then $\dist_{H^*}(f^*(u), f^*(v)) \geq \|h(u) - h(v)\|_\infty > 0$, as $h$ is $R$-locally injective. On the other hand, if $\dist_G(u,v) \geq R$, then $\dist_{H^*}(f^*(u), f^*(v)) \geq \dist_H(f(u),f(v)) \geq \dist_G(u,v)^{1-\epsilon} > 0$ by \eqref{eq:far}. Therefore, $f^*$ is an embedding of $G$ into $H^*$. To finish the proof of Theorem~\ref{theo:coarse}, it remains to observe that $\rho \leq \rho(\gr{d}) = \Theta(d)$, and hence $n + 4d = O_\epsilon(d)$ and $R$ is bounded above by a function of $d$ and $\epsilon$.

    \section{Proof of Lemma~\ref{lemma:key}}\label{sec:asi}

    \subsection{Borel asymptotic dimension and asymptotic separation index}

    One of the fundamental tools in large-scale geometry is the concept of asymptotic dimension of a metric space, introduced by Gromov \cite[\S1.E]{gromov1993asymptotic}. In their highly influential paper \cite{Dimension}, Conley, Jackson, Marks, Seward, and Tucker-Drob developed a Borel version of this notion. They also defined a closely related parameter, the asymptotic separation index, which turned out to be extremely useful in Borel combinatorics \cite{BWKonig,Dimension,ASIalgorithms,FelixFinDim,BernshteynWeilacher}.
    
    For a graph $G$ and an integer $R \in \N$, $G^R$ denotes the graph with $V(G^R) \defeq V(G)$ in which two vertices $u$, $v$ are adjacent if and only if $0 < \dist_G(u,v) \leq R$. If $G$ is a locally finite Borel graph, then $G^R$ is also locally finite and Borel \cite[Corollary 5.2]{Pikh_survey}. We say that a subset $U \subseteq V(G)$ is \emphd{$G$-finite} if every component of the induced subgraph $G[U]$ is finite, and \emphd{$G$-bounded} if there is $D \in \N$ such that every component of $G[U]$ has $\dist_G$-diameter at most $D$.

    \begin{defn}[Borel asymptotic dimension and asymptotic separation index]
        Let $G$ be a locally finite Borel graph. We define parameters 
        $\asi(G)$ and $\asdimB(G)$ as follows.
        \begin{itemize}
             \item The \emphd{asymptotic separation index} of $G$, in symbols $\asi(G)$, is the smallest $s \in \N$ such that for every $R \in \N$, there exists a cover $V(G) = U_0 \cup \ldots \cup U_s$ of $V(G)$ by $G^R$-finite Borel subsets. If no such $s \in \N$ exists, we set $\asi(G) \defeq \infty$.

             \item The \emphd{Borel asymptotic dimension} of $G$, in symbols $\asdimB(G)$, is the smallest $s \in \N$ such that for every $R \in \N$, there exists a cover $V(G) = U_0 \cup \ldots \cup U_s$ of $V(G)$ by $G^R$-bounded Borel subsets. If no such $s \in \N$ exists, we set $\asdimB(G) \defeq \infty$.
            



        \end{itemize}
    \end{defn}

    It is clear that $\asdimB(G) \geq \asi(G)$. In \cite{Dimension}, Conley \emph{et al.}~gave numerous examples of Borel graphs with finite Borel asymptotic dimension. For instance, they showed that Schreier graphs of Borel actions of various finitely generated groups (polycyclic groups, the lamplighter group $\Z_2 \wr \Z$, and the Baumslag--Solitar group $BS(1,2)$, to name a few) have finite Borel asymptotic dimension. They also established the following relationship between $\asdimB(G)$ and $\asi(G)$:

    \begin{theo}[{Conley--Jackson--Marks--Seward--Tucker-Drob \cite[Theorem 4.8(a)]{Dimension}}]\label{theo:asdim_to_asi}
        If $G$ is a locally finite Borel graph with $\asdimB(G) < \infty$, then $\asi(G) \leq 1$.
    \end{theo}

    One of the results of \cite{BernshteynYu} is a bound on the Borel asymptotic dimension for Borel graphs of polynomial growth (which is analogous to a bound on the ordinary asymptotic dimension for graphs of polynomial growth due to Papasoglu \cite{papasoglu2021polynomial}):

    \begin{theo}[{AB--JY \cite[Theorem 1.25]{BernshteynYu}}]\label{theo:polygrowth_asdim}
        Let $G$ be a Borel graph of polynomial growth. Then the Borel asymptotic dimension of $G$ is finite. Moreover, if for all sufficiently large $r \in \N$, every $r$-ball in $G$ contains at most $(r+1)^\rho$ vertices, then $\asdimB(G) \leq \rho$; in particular, $\asdimB(G) \leq \rho(G)$.
    \end{theo}

    In our proof of Lemma~\ref{lemma:key}, we use the following consequence of Theorems~\ref{theo:asdim_to_asi} and \ref{theo:polygrowth_asdim}:

    \begin{corl}\label{corl:asi}
        If $G$ is a Borel graph of polynomial growth, then $\asi(G) \leq 1$.
    \end{corl}

    \subsection{Constructing locally injective Lipschitz functions}

    An indispensable tool when working with locally finite (or, more generally, locally countable) graphs is the Luzin--Novikov theorem \cite{Luzin,Novikov}:

    \begin{theo}[{Luzin--Novikov \cite[Theorem 18.10]{KechrisDST}}]\label{theo:LN}
        Let $X$ and $Y$ be standard Borel spaces and let $A \subseteq X \times Y$ be a Borel subset. If for all $x \in X$, the set $A_x \defeq \set{y \in Y \,:\, (x,y) \in A}$ is countable, then there exists a sequence $(f_n)_{n \in \N}$ of Borel partial functions $f_n \colon X \pto Y$ defined on Borel subsets of $X$ such that for all $x \in X$, $A_x = \set{f_n(x) \,:\, n \in \N,\ x \in \dom(f_n)}$.
    \end{theo}

    A consequence of the Luzin--Novikov theorem is, roughly, that constructions that only involve quantifiers ranging over countable sets produce Borel results. For example, it implies that for a locally finite Borel graph $G$, the functions $\deg_G \colon V(G) \to \N$ and $\dist_G \colon V(G)^2 \to \N \cup \set{\infty}$ are Borel. For a detailed discussion and further examples of using Theorem~\ref{theo:LN} in Borel graph combinatorics, see \cites[\S3.1]{BernshteynYu}[\S5]{Pikh_survey}. Routine arguments involving the Luzin--Novikov theorem similar to the ones described in the cited references will be employed in the sequel without mention. Another useful consequence of the Luzin--Novikov theorem is the following \emphd{uniformization} result: In the setting of Theorem~\ref{theo:LN}, if each set $A_x$ is nonempty, then the mapping $f \colon X \to Y$ given by
        \[
            f(x) \,\defeq\, f_{n}(x), \text{ where $n \in \N$ is minimum such that $x \in \dom(f_n)$},
        \]
    is a Borel function such that $f(x) \in A_x$ for all $x \in X$, i.e., $f$ selects one element from each $A_x$.
    
    Before we commence the proof of Lemma~\ref{lemma:key}, we need to establish a few auxiliary facts. 
    Let $\mathcal{I}_R$ denote the interval $\set{0,1,\ldots, R}$ of length $R$ in $\Z$. We view the $d$-dimensional box $\mathcal{I}_R^d$ as a subset of $\Z^d$ equipped with the $\infty$-metric.

    \begin{lemma}\label{lemma:fold}
        For any $d$, $R \in \N$, there exists an $R$-locally injective $1$-Lipschitz function $f \colon \Z^d \to \mathcal{I}_R^{2d}$.
    \end{lemma}
    \begin{scproof}
        If $f \colon \Z \to \mathcal{I}_R^2$ is an $R$-locally injective $1$-Lipschitz function, then so is the map $f^{d} \colon \Z^d \to \mathcal{I}_R^{2d}$ given by $f^{d}(k_1, \ldots, k_d) \defeq (f(k_1), \ldots, f(k_d))$. (Here it is crucial that we are using the $\infty$-metric.) Therefore, it is enough to prove the lemma for $d = 1$. To this end, note that there exists a $4R$-cycle in the graph $\gr{2}$ that is entirely contained in $\mathcal{I}_R^2$, namely the $R$-by-$R$ square with corners $(0,0)$, $(R,0)$, $(R,R)$, and $(0,R)$. Let $v_0$, $v_1$, \ldots, $v_{4R-1} \in \mathcal{I}_R^2$ be the vertices of this cycle in their cyclic order and define $f \colon \Z \to \mathcal{I}_R^2$ by $f(k) \defeq v_{k \, \mathrm{mod} \, 4R}$. The function $f$ is $1$-Lipschitz and $R$-locally injective (indeed, it is $(4R-1)$-locally injective), so we are done.
    \end{scproof}


    
    \begin{corl}\label{corl:fold}
        If $G$ is a graph with $\emb(G) = d < \infty$, then for any $R \in \N$ and every finite subset $S \subseteq V(G)$, there exists a function $f \colon S \to \mathcal{I}_R^{2d}$ that is $R$-locally injective and $1$-Lipschitz with respect to the metric $\dist_G$ restricted to $S$.
    \end{corl}
    \begin{scproof}
        Let $D$ be the maximum of the finite distances between the vertices in $S$ and let $F$ be the subgraph of $G$ induced by the vertices at distance at most $D$ from $S$ (here all distances are in the metric $\dist_G$). Then $\dist_{F}(u,v) = \dist_G(u,v)$ for all $u$, $v \in S$. Since $G$ is locally finite, $F$ is a finite subgraph of $G$, so, by the definition of $\emb(G)$, there exists an embedding from $F$ into $\gr{d}$. Composing this embedding with the $R$-locally injective $1$-Lipschitz function $\Z^d \to \mathcal{I}_R^{2d}$ given by Lemma~\ref{lemma:fold} yields a mapping $f \colon V(F) \to \mathcal{I}_R^{2d}$ that is $R$-locally injective and $1$-Lipschitz with respect to $\dist_{F}$. As $\dist_{F}$ agrees with $\dist_G$ on $S$, taking the restriction of $f$ to $S$ finishes the proof.
    \end{scproof}

    The next lemma allows us to extend a $1$-Lipschitz function with codomain $\Z^d$ from a subset of $V(G)$ to the entire vertex set. This problem is considerably simplified by the fact that we equip $\Z^d$ with the $\infty$-metric. For example, the same question for the $1$-metric was investigated by Chandgotia, Pak, and Tassy \cite{LipExt} and is significantly more intricate.

    \begin{lemma}\label{lemma:extend}
        Let $G$ be a locally finite Borel graph. Suppose that $A \subseteq V(G)$ is a Borel subset and $f \colon A \to \Z^d$ is a Borel function that is $1$-Lipschitz with respect to the metric $\dist_G$ restricted to $A$. Then $f$ can be extended to a $1$-Lipschitz Borel function $f^* \colon V(G) \to \Z^d$.
    \end{lemma}
    \begin{scproof}
        Since we are using the $\infty$-metric on $\Z^d$, a function $f$ from a metric space $(X, d_X)$ to $\Z^d$ is $1$-Lipschitz if and only if for all $1 \leq i \leq d$, the map $f_i \colon X \to \Z$ sending each $x$ to the $i$-th coordinate of $f(x)$ is $1$-Lipschitz. Hence, it is enough to prove the lemma in the case $d = 1$. So, let $A \subseteq V(G)$ be a Borel subset and let $f \colon A \to \Z$ be a $1$-Lipschitz (with respect to $\dist_G$) Borel function. Let $[A]_G$ be the set of all vertices $u \in V(G)$ such that $\dist_G(u,A) < \infty$. Since $G$ is locally finite, $[A]_G$ is a Borel subset of $V(G)$ by the Luzin--Novikov theorem (Theorem~\ref{theo:LN}). By definition, $G$ has no edges joining $[A]_G$ to $V(G) \setminus [A]_G$. Now, for each $u \in V(G)$, we define
        \[
            f^*(u) \,\defeq\, \begin{cases}
                \min_{a \in A} (f(a) + \dist_G(u,a)) &\text{if $u \in [A]_G$},\\
                0 &\text{if $u \notin [A]_G$}.
            \end{cases}
        \]
        For $u \in [A]_G$, there are only countably many vertices $a \in A$ with $\dist_G(u,a) < \infty$, so the minimum in the definition of $f^*(u)$ is taken over a countable set. By the Luzin--Novikov theorem again, $f^*$ is Borel. It is also straightforward to verify that $f^*$ is $1$-Lipschitz and extends $f$. 
    \end{scproof}

    We now have all the ingredients required to prove Lemma~\ref{lemma:key}, which we restate here for the reader's convenience:

    \begin{lemmacopy}{lemma:key}
        Let $G$ be a Borel graph of polynomial growth with $\emb(G) = d$. Then, for any $R > 0$, there exists an $R$-locally injective $1$-Lipschitz Borel function $f \colon V(G) \to \Z^{4d}$.
    \end{lemmacopy}
    \begin{scproof}
        By Corollary~\ref{corl:asi}, $\asi(G) \leq 1$. Therefore, there exists a cover $V(G) = U_0 \cup U_1$ of $V(G)$ by a pair of $G^{3R}$-finite Borel sets. For $i \in \set{0,1}$, let
        \[
            V_i \,\defeq\, \set{v \in V(G) \,:\, \dist_G(v, U_i) \leq R}.
        \]
        Note that if $v$, $v' \in V_i$ are two vertices with $\dist_G(v,v') \leq R$ and $u$, $u' \in U_i$ are such that $\dist_G(v,u)$, $\dist_G(v',u') \leq R$, then $\dist_G(u,u') \leq 3R$ and hence $u$ and $u'$ are in the same component of $(G^{3R})[U_i]$. Since $U_i$ is $G^{3R}$-finite and $G$ is a locally finite graph, it follows that $V_i$ is $G^R$-finite.

        Consider any connected component $H$ of the graph $(G^R)[V_i]$. By the above discussion, $H$ is finite, so, by Corollary~\ref{corl:fold}, there exists a function $f_H \colon V(H) \to \mathcal{I}_R^{2d}$ that is $R$-locally injective and $1$-Lipschitz with respect to the metric $\dist_G$ restricted to $V(H)$. Furthermore, there are only finitely many choices for such $f_H$, so by the Luzin--Novikov uniformization, the mapping $H \mapsto f_H$ can be arranged to be Borel \ep{this is a special case of the general fact that on a Borel graph with finite components, a local labelling problem that has any solution has a Borel solution as well; see \cites[\S2.2]{BernshteynWeilacher}[\S5.3]{Pikh_survey}}. Now we define a Borel function $f_i \colon V_i \to \mathcal{I}_R^{2d}$ by
        \[
            f_i(v) \,\defeq\, f_H(v), \text{ where $H$ is the component of $(G^R)[V_i]$ containing $v$}.
        \]
        We claim that the map $f_i$ is $1$-Lipschitz with respect to the metric $\dist_G$ restricted to $V_i$. Indeed, if $v$, $v' \in V_i$ are in the same component $H$ of $(G^R)[V_i]$, then
        \[
            \|f_i(v) - f_i(v')\|_\infty \,=\, \|f_H(v) - f_H(v')\|_\infty \,\leq\, \dist_G(v,v')
        \]
        by the choice of $f_H$. On the other hand, if $v$ and $v'$ are in different components of $(G^R)[V_i]$, then $\dist_G(v,v') > R \geq \|f_i(v) - f_i(v')\|_\infty$, because the diameter of $\mathcal{I}_R^{2d}$ in the $\infty$-metric is $R$. Therefore, by Lemma~\ref{lemma:extend}, we can extend $f_i$ to a Borel $1$-Lipschitz function $f^*_i \colon V(G) \to \Z^{2d}$.

        Finally, we define $f \colon V(G) \to \Z^{4d}$ by $f(v) \defeq (f^*_0(v), f^*_1(v))$ for all $v \in V(G)$. Since $f^*_0$ and $f^*_1$ are $1$-Lipschitz, so is $f$. Now take any $u$, $v \in V(G)$ with $0 < \dist_G(u,v) \leq R$. Since $U_0 \cup U_1 = V(G)$, we have $u \in U_i$ for some $i \in \set{0,1}$. Then $u$ and $v$ both belong to $V_i$ and, moreover, they are in the same component, say $H$, of $(G^R)[V_i]$. Since $f_H$ is $R$-locally injective, it follows that $f_i^*(u) = f_H(u) \neq f_H(v) = f_i^*(v)$, and thus $f(u) \neq f(v)$. Therefore, $f$ is $R$-locally injective, and we are done.
    \end{scproof}

    \section{Further remarks}\label{sec:further}

    Here we record some further facts about Borel embeddings into Borel $\gr{n}$-graphs, which help clarify the relationship between different natural approaches to defining Borel embedding dimension.

    The following definition will be useful: Given a locally finite Borel graph $G$, we let $\sim_G$ be the \emphd{connectedness relation} of $G$, i.e., the equivalence relation on $V(G)$ whose classes are the connected components of $G$; in other words, $u \sim_G v$ if and only if $\dist_G(u,v) < \infty$. Note that the $\sim_G$-equivalence classes are countable. As a subset of $V(G)^2$, $\sim_G$ is Borel \cite[Corollary 5.3]{Pikh_survey}.

    We say that an embedding $f \colon V(G) \to V(H)$ of a graph $G$ into a graph $H$ is \emphd{strong} if vertices in different components of $G$ are mapped to distinct components of $H$. When $H$ is a Borel $\gr{d}$-graph, this means that every component of $G$ is embedded into its own copy of $\gr{d}$. The following proposition shows that using strong embeddings would make no difference in Definition~\ref{defn:Borel_dim}:

    \begin{prop}\label{prop:strong}
        Let $G$ be a Borel graph. The following statements are equivalent:
        \begin{enumerate}[label=\ep{\normalfont\roman*}]
            \item\label{item:emb} $\embB(G) \leq d$,
            \item\label{item:strong_emb} $G$ admits a Borel strong embedding into a Borel $\gr{d}$-graph.
        \end{enumerate}
    \end{prop}
    \begin{scproof}
        The implication \ref{item:strong_emb} $\Longrightarrow$ \ref{item:emb} is clear. To prove \ref{item:emb} $\Longrightarrow$ \ref{item:strong_emb}, suppose that $\embB(G) \leq d$ and let $f \colon V(G) \to V(H)$ be a Borel embedding of $G$ into a Borel $\gr{d}$-graph $H$ generated by a free Borel action $\Z^d \acts V(H)$. Since $f$ is an embedding, if $u \sim_G v$, then $f(u) \sim_H f(v)$, and hence there exists a unique group element $\updelta(u,v) \in \Z^d$ such that $\updelta(u,v) \bm{+} f(u) = f(v)$. Note that the map $\updelta \colon {\sim_G} \to \Z^d$ has the following properties:
        \begin{itemize}
            \item $\updelta(u,u) = 0$ for all $u \in V(G)$,
            \item $\updelta(u,v) + \updelta(v,w) = \updelta(u,w)$ for all $u \sim_G v \sim_G w$,
            \item the relation $E_\updelta \defeq \set{(u,v) \in V(G)^2 \,:\, \updelta(u,v) = 0}$ is trivial.
        \end{itemize}
        By \cite[Theorem 4.5]{BernshteynYu}, this implies that there exist a free Borel action $\Z^d \acts Y$ on some standard Borel space $Y$ and a Borel map $f^* \colon V(G) \to Y$ such that for all $u$, $v \in V(G)$,
        \begin{itemize}
            \item if $u \sim_G v$, then $\updelta(u,v) \bm{+} f^*(u) = f^*(v)$, and
            \item if $u \not\sim_G v$, then $f^*(u)$ and $f^*(v)$ are in different orbits of the action $\Z^d \acts Y$.
        \end{itemize}
        Let $H^* \defeq \Sch(Y, S_{d,\infty})$ be the Borel $\Z^d_\infty$-graph corresponding to the action $\Z^d \acts Y$. We claim that $f^*$ is a Borel strong embedding of $G$ into $H^*$. Indeed, take any $u$, $v \in V(G)$. If $u \not\sim_G v$, then $f^*(u)$ and $f^*(v)$ are in different orbits of the action $\Z^n \acts Y$, hence $\dist_{H^*}(f^*(u),f^*(v)) = \infty$. On the other hand, if $u \sim_G v$, then $\updelta(u,v) \bm{+} f^*(u) = f^*(v)$, which implies that \[\dist_{H^*}(f^*(u),f^*(v)) \,=\, \|\updelta(u,v)\|_\infty \,=\, \dist_H(f(u), f(v)).\] Since $f$ is an embedding, it follows that $f^*$ is $1$-Lipschitz and injective, as desired.
    \end{scproof}
    
    Certain Borel $\gr{n}$-graphs are ``universal'' from the point of view of Borel embeddings. Given a standard Borel space $X$, we let $X^{\Z^n}$ be the product space of all functions $x \colon \Z^n \to X$ and define the \emphd{\ep{Bernoulli} shift action} $\Z^n \acts X^{\Z^n}$ via the formula
    \[
        (k \bm{+} x)(m) \,\defeq\, x(m + k) \quad \text{for all $x \in X^{\Z^n}$ and $k$, $m \in \Z^n$}.
    \]
    The \emphd{free part} of $X^{\Z^n}$, denoted by $\Free(X^{\Z^n})$, is the set of all points $x \in X^{\Z^n}$ whose stabilizer under the shift action is trivial. By definition, $\Free(X^{\Z^n})$ is the largest shift-invariant subspace of $X^{\Z^n}$ on which the shift action is free. We denote the Borel $\Z^n_\infty$-graph generated by 
    the restricted shift action $\Z^n \acts \Free(X^{\Z^n})$ by $\mathsf{ShiftGrid}^n_\infty(X)$.

    \begin{prop}\label{prop:universal}
        Let $H$ be a Borel $\gr{n}$-graph.

        \begin{enumerate}[label=\ep{\normalfont\arabic*}]
            \item\label{item:R} $H$ admits a Borel distance-preserving embedding into $\shgr_\infty^{n}(\R)$.

            \item\label{item:two} $H$ admits a Borel distance-preserving embedding into $\shgr_\infty^{n+1}(\set{0,1})$.
        \end{enumerate}
    \end{prop}
    \begin{scproof}
        For part \ref{item:R}, thanks to the Borel isomorphism theorem \cite[Theorem 15.6]{KechrisDST}, we may assume that $V(H) = \R$. Let $\Z^n \acts \R$ be the free Borel action that generates $H$. The desired embedding is given by the map $f \colon \R \to \R^{\Z^n}$ defined as follows:
        \[
            (f(x))(k) \,\defeq\, k \bm{+} x \quad \text{for all $x \in \R$ and $k \in \Z^n$}.
        \]
        Part \ref{item:two} is verified in \cite[proof of Lemma 4.4]{BernshteynYu}.
    \end{scproof}

    \subsection*{Acknowledgment} We are grateful to the anonymous referee for their helpful feedback.

	\printbibliography
    
\end{document}